\documentclass[12pt,dvips]{amsart}
\usepackage{euler, amsfonts, amssymb, latexsym, epsfig,epic}

\newcommand\lie[1]{{\mathfrak #1}}

\newcommand\Tr{{\rm Tr\,}}

\newcommand\Spec{{\rm Spec\,}}
\newcommand\Proj{{\rm Proj\,}}
\newcommand\Sym{{\rm Sym\,}}
\newcommand\iso{{\cong}}
\newcommand\tensor{{\otimes}}

\newtheorem{Theorem}{Theorem} 
\newtheorem{Proposition}{Proposition}

\newtheorem*{Corollary}{Corollary}
 
\newtheorem*{Proposition*}{Proposition} 
\newtheorem*{Theorem*}{Theorem}
\theoremstyle{remark}
\newtheorem{Example}{Example}

\newcommand\onto{\mathop{\to\!\!\!\!\to}}
\newcommand\into{\operatorname*{\hookrightarrow}}

\newcommand\CP{{\mathbb C \mathbb P}}
\newcommand\PP{{\mathbb P}}

\newcommand\reals{{\mathbb R}}
\newcommand\complexes{{\mathbb C}}
\newcommand\integers{{\mathbb Z}}

\newcommand\naturals{{\mathbb N}}

\theoremstyle{plain}

\newcommand\dfn{\bf} 

\newcommand\GLn{{{GL_n(\complexes)}}}

\begin{document}
\pagestyle{plain}

\title{The symplectic and algebraic geometry of Horn's problem}
\author{Allen Knutson}
\email{allenk@math.berkeley.edu}
\date{\today}

\maketitle

\begin{abstract}
  One version of Horn's problem asks for which $\lambda,\mu,\nu$ does
  $H_\lambda + H_\mu + H_\nu = 0$ have solutions, where
  $H_{\lambda,\mu,\nu}$ are Hermitian matrices with spectra
  $\lambda,\mu,\nu$. This turns out to be a {\em moment map} condition
  in Hamiltonian geometry. Many of the results around Horn's problem
  proven with great effort ``by hand'' are in fact simple consequences
  of the modern machinery of symplectic geometry, and the subtler ones
  provable via the connection to geometric invariant theory. We give
  an overview of this theory (which was not available to Horn),
  including all definitions, and how it can be used in linear algebra.
\end{abstract}

\section{Introduction}

This is an expository paper on the symplectic and algebraic geometry
implicit in {\em Horn's problem}, which asks the possible spectra of a
sum of two Hermitian matrices each with known spectrum.

The connection with symplectic geometry is very straightforward: 
the map $(H_\lambda,H_\mu) \mapsto H_\lambda+H_\mu$ that takes a pair
of Hermitian matrices with known spectra to their sum is a
{\em moment map} for the diagonal conjugation action of $U(n)$ on
a certain {\em symplectic manifold} (definitions to follow). 
This is a very restrictive 
property of maps, and many things can be proved about them.
The proofs, at heart, are not really any different than the
techniques Horn himself used to study this map. Nonetheless 
the framework is worth understanding in order to recognize 
what other linear algebra problems are likely to have answers 
as nice as the ones to Horn's problem. In particular the Schur-Horn
theorem follows very easily from the general theorems in this area
(and was a primary inspiration for them).

Some of the more esoteric connections -- to algebraic geometry and
representation theory -- can also be seen in this context, via
the Kirwan/Ness theorem (which we will also state). 
Again, the basic techniques used are the same,
but in the Kirwan/Ness theorem one sees these techniques pushed to
prove the statements in what appears to be their proper generality.

Along the way we explain the relation between Hermitian matrices,
flag manifolds, and the Borel-Weil-Bott-Kostant theorem.

\section{The Schur-Horn theorem, Horn's problem, and Hamiltonian manifolds}

Let $\lambda = (\lambda_1 \geq \lambda_2 \geq \ldots \geq \lambda_n)$
be a weakly decreasing list of real numbers, which we'll use to encode
the eigenvalue spectrum of a Hermitian matrix.
The Schur-Horn theorem states the following:

\begin{Theorem*}
  Let ${\mathcal O}_\lambda$ be the space of Hermitian matrices with
  spectrum $\lambda$. Let $\Phi : {\mathcal O}_\lambda \to \reals^n$ take
  a matrix to its diagonal entries. Then the image of $\Phi$ is a 
  convex polytope, whose vertices are the $n!$ permutations of $\lambda$.
\end{Theorem*}

This theorem was interpreted by Kostant in 1970 as the
$U(n)$ case of a theorem for arbitrary compact Lie groups, 
leading the way to a much wider generalization found in 1982 by Atiyah and 
independently by Guillemin and Sternberg:

\begin{Theorem*}
  Let $M$ be a compact connected symplectic manifold, with an action
  of a torus $T$. Let $\Phi : M \to \lie{t}^*$ be a moment map for
  this action. Then the image of $\Phi$ is a convex polytope, the
  convex hull of the images of the $T$-fixed points on $M$.
\end{Theorem*}

Of course, to see how to cast the Schur-Horn theorem in this formulation,
we'll need to define ``symplectic manifold'' and ``moment map''.
We will make no attempt to be encyclop\ae dic in our references and
instead direct the reader to 
\begin{itemize}
\item \cite{GLS} -- for all matters symplectic or convex
\item \cite{F} -- for the algebraic geometry of flag manifolds
\item \cite{MFK} -- for geometric invariant theory
\item chapter 8 of \cite{MFK} -- for the Kirwan/Ness theorem
\end{itemize}
and references therein.

\subsection{Symplectic manifolds.} Let $M$ be a manifold, and $\omega$
an {\em anti}symmetric inner product on the tangent spaces to $M$,
sort of a skew Riemannian metric.
Since the inner product of any vector with itself is zero, we can't talk
about positive definiteness, and so we instead ask for nondegeneracy --
that for any tangent vector $\vec v_1$, there exists another vector $\vec v_2$
with $\omega(\vec v_1,\vec v_2)$ nonzero. 
Surprisingly, this forces $M$ to be even dimensional.

There is a standard example: $\reals^{2d}$ with basis 
$\{\vec x_1,\ldots,\vec x_d,\vec y_1,\ldots,\vec y_d\}$, where 
$\omega(\vec x_i,\vec x_j)=\omega(\vec y_i,\vec y_j)=0$, but
$\omega(\vec x_i,\vec y_j) = - \omega(\vec y_j,\vec x_i) = 
\delta_{ij}$ (Kronecker delta).

If $(M,\omega)$ is locally isomorphic to $\reals^{2d}$ with its
standard $\omega$, we say that $M$ is a {\dfn symplectic manifold},
and call $\omega$ its {\dfn symplectic form}.\footnote{%
The usual definition is in terms of de Rham forms, and has this
formulation as Darboux's theorem.} This is roughly analogous to
studying Riemannian manifolds that are locally Euclidean.

As with a Riemannian metric, we can talk about the {\dfn symplectic gradient}
$X_f$ of a function $f$, also called the Hamiltonian vector field $X_f$ 
associated to a Hamiltonian $f$. It is defined uniquely by the equation
$$ D_{\vec v(m)} f = \omega(\vec v(m), X_f(m)) $$
where $\vec v(m), X_f(m)$ are tangent vectors to the point $m\in M$, and
$D_{\vec v(m)} f$ is the directional derivative of $f$ in the
direction $\vec v(m)$. That this defines $X_f$ uniquely follows from 
$\omega$'s nondegeneracy.

It is easy to show that symplectic gradients have an immense advantage
over Riemannian gradients: the derivative of $\omega$ along $X_f$ vanishes.
In integral form this says that the time $t$ flow map from $M$ to $M$,
given by following $X_f$ for time $t$, takes $\omega$ to itself.

\begin{figure}[htbp]
  \begin{center}
    \epsfig{file=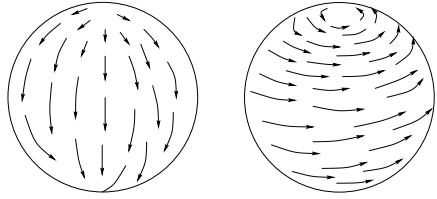,height=1in}
    \caption{The Riemannian gradient vs. the symplectic gradient of the
      height function on $S^2$.}
    \label{fig:example1}
  \end{center}
\end{figure}

\begin{Example}
  Let $M=S^2$, and $\omega$ be the area form, taking a pair of vectors 
  to the oriented area of the parallelogram they define. Let $f$ be
  the height function on $M$, normalized to take the North pole to $1$
  and the South pole to $-1$. Then the Riemannian gradient points down
  everywhere along longitude lines, as in figure \ref{fig:example1}, 
  whereas the symplectic gradient points sideways everywhere along
  latitude lines -- and generates rotations, which are of course
  area-preserving.
\end{Example}

This gives us a handy way to describe (certain) flows on a symplectic
manifold -- in terms of functions, which opens up the world 
of e.g. Morse theory. One example of this connection: 
a critical point of a function $f$ is one where all 
the directional derivatives are zero, which (by $\omega$'s nondegeneracy)
is equivalent to the vanishing there of $X_f$. So the critical
points of $f$ are exactly the fixed points of the flow generated by
its symplectic gradient $X_f$.

\subsection{Moment maps.}
Now let $M$ be a symplectic manifold with symplectic form $\omega$,
and an action of a connected Lie group $K$;  for us $K$ will always
be the circle group $S^1$, a product $T$ of circle groups, or the
unitary group $U(n)$.

If we assume $K$ acts smoothly on $M$, 
the elements of $K$ nearby the identity give
us diffeomorphisms of $M$ very close to the identity; differentiating
this picture, we find that each tangent vector to the identity of $K$
gives a vector field on $M$. Denote the tangent space to the identity,
also known as the Lie algebra of $K$, by $\lie{k}$, and its dual
by $\lie{k}^*$. (This funny letter is a Fraktur $k$.)

Since $K$ acts on itself by conjugation, fixing the identity, it
acts on $\lie{k}$; this is called the {\dfn adjoint representation},
and the induced action on $\lie{k}^*$ the {\dfn coadjoint representation}.

We say that a map $\Phi : M \to \lie{k}^*$ is a {\dfn moment map} for
the action of $K$ on $M$ if
\begin{enumerate}
\item $\Phi$ is equivariant, i.e. $  \forall k\in K,m\in M,
  \Phi(k\cdot m) = k \Phi(m) k^{-1}$
\item for each $\vec k \in \lie{k}$, the vector field induced on $M$
  by $\vec k$ equals the symplectic gradient of $\langle \vec k, \Phi\rangle$.
  (We are using here the natural pairing between $\lie{k}$ and $\lie{k}^*$.)
\end{enumerate}

Note that for the action to have a moment map at all, the action of
$K$ on $M$ must preserve the symplectic structure. This is almost
sufficient (but not quite); in any case all the actions considered in
this paper will have moment maps. If the action has a moment map
it is said to be {\dfn Hamiltonian} and $M$ is called a
{\dfn Hamiltonian $K$-manifold}.

Since $K$ is connected, we can recover its action from the action of
the Lie algebra, which in turn we can get from the moment map;
in particular the moment map uniquely determines the action.
(The action does not determine the moment map uniquely even when
one exists, but only up to certain translations.)
And while any map $\Phi:M\to \lie{k}^*$ will give us associated
vector fields $X_{\langle \vec k,\Phi\rangle}$, very few maps
will give actions of $K$; moment maps are very special.

We need four particularly important facts about moment maps.
(All follow easily from the definition.)
\begin{enumerate}
\item Let $K$ act on $M$ with $\Phi_K : M \to \lie{k}^*$ a
  moment map for the action on $M$, and $\rho: H \to K$ be a Lie group
  homomorphism, making $H$ act on $M$ too.  
  There is a corresponding map of Lie algebras $d\rho: \lie{h} \to \lie{k}$,
  and a dual map $(d\rho)^*: \lie{k}^* \to \lie{h}^*$. Then the action of $H$
  is also Hamiltonian, with moment map given by composing $\Phi_K$
  and $(d\rho)^*$. The case of greatest interest is when $H$ is a
  subgroup of $K$, and $\rho$ is inclusion.
\item Let $M$ be a {\dfn coadjoint orbit}, i.e. an orbit of $K$'s 
  action on $\lie{k}^*$. Then there exists uniquely a symplectic structure
  on $M$ such that the inclusion map $M \into \lie{k}^*$ is a moment map.
  (This was called the Kirillov-Kostant-Souriau symplectic structure until
  Alan Weinstein found it in 19th-century notebooks of Lie.)
\item Let $M$ be a Hamiltonian $K$-manifold, and $N$ 
  a Hamiltonian $H$-manifold. Then the natural action of $K\times H$
  on $M\times N$ is Hamiltonian, with moment map the direct sum of
  the two individual moment maps.
\item Let $M$ be a Hamiltonian $K$-manifold, and $N$ a submanifold of
  $M$ invariant under $K$, such that the restriction of the symplectic
  form on $M$ is nondegenerate on $N$ (making $N$ a symplectic manifold
  in a natural way). Then the action of $K$ on $N$ is also Hamiltonian,
  with moment map the composition of inclusion with $M$'s own moment map.
\end{enumerate}

Since the image of a moment map is $K$-invariant, it must be a union
of coadjoint orbits; in this way one can regard the individual coadjoint
orbits as sort of ``minimal'' Hamiltonian $K$-manifolds. (This is only
interesting for nonabelian groups $K$, insofar as the coadjoint action 
is trivial for abelian groups, with the coadjoint orbits just points.)

\begin{Example}
  Let $M=\reals^2$ with the standard symplectic form, and $K=S^1$
  acting on $M$ by rotation. Then we can identify $K$'s tangent space
  with $\reals$, such that the moment map is $\Phi(\vec v) = |\vec v|^2$.
  The image of $\Phi$ is the positive half-line.
\end{Example}

\begin{Example}[using property 3 above of moment maps]
  Let $M=\reals^{2d}$ with the standard symplectic form, and $K=(S^1)^d$,
  each circle acting on a pair of the coordinates -- really the
  $n$th power of the previous example.
  Then we can identify $K$'s tangent space with $\reals^d$.
  The image of $\Phi$ is the positive ``orthant'' (the $d$-dimensional
  generalization of the first quadrant in $\reals^2$, octant in $\reals^3$,
  etc.).
\end{Example}

This second example gives a hint of the source of convex polytopes
in the Atiyah/Guillemin-Sternberg theorem stated above; since we've
defined a symplectic manifold as one that looks locally like $\reals^{2d}$
with the standard form, we just need to prove a ``$T$-equivariant version''
of that local normal form, in order to know that at least locally the
image of the moment map is a polytope. The actual proof of A/G-S 
requires some more ingredients (not explored here) 
to show the global statement.

\begin{Example}
  Let $M=\CP^n=\{ [z_0,z_1,\ldots,z_n] : 
  z_i \hbox{ not all } 0\}/\complexes^\times$, complex projective $n$-space.
  This has a natural symplectic form (the Fubini-Study form)
  we do not pause to write down, and a Hamiltonian action of
  $T^{n+1}$, whose Lie algebra we can identify with $\reals^{n+1}$.
  One moment map is then
  $$  \Phi([z_0,z_1,\ldots,z_n]) = \bigg( \frac{|z_0|^2}{\sum_i |z_i|^2},
  \frac{|z_1|^2}{\sum_i |z_i|^2},\ldots,
  \frac{|z_n|^2}{\sum_i |z_i|^2} \bigg). $$
  The individual coordinates can only vary between $0$ and $1$, and their
  sum is automatically $1$; it is easy to check that the image is the
  whole simplex. In fact this is really the $\lambda=(1,0,0,\ldots,0)$
  case of the Schur-Horn result.
\end{Example}

\subsection{The relation with the Schur-Horn theorem.}
Let $K$ be the unitary group $U(n)$, and $T$ the subgroup of diagonal
matrices, a ``maximal torus'' of $U(n)$. The Lie algebra $\lie{u}(n)$ is
the space of skew-Hermitian matrices.\footnote{%
A unitary matrix has $U U^* = 1$. A unitary matrix near the identity,
$1 + \epsilon S$, therefore has $(1 + \epsilon S)(1 + \epsilon S)^* =
1 + \epsilon(S+S^*) + O(\epsilon^2) = 1$. Differentiating, we get
$S = -S^*$.}
We can identify the Hermitian matrices with $\lie{u}(n)^*$ 
using the trace form, $H \mapsto \Tr (i H \cdot)$. This identification
is $U(n)$-equivariant, intertwining the conjugation action with
the coadjoint action; in particular it takes orbits of Hermitian
matrices to coadjoint orbits.

The upshot is that we can think of coadjoint orbits of $U(n)$ as
orbits of $U(n)$ acting on the space of Hermitian matrices by conjugation, 
which of course are the isospectral sets ${\mathcal O}_\lambda$
required for the Schur-Horn theorem.

From this point, we just have to apply properties 1 and 2 of moment maps;
${\mathcal O}_\lambda$ is a symplectic manifold, and the action of 
the diagonal subgroup $T$ is Hamiltonian with moment map 
${\mathcal O}_\lambda \into \lie{u}(n)^* \onto \lie{t}^*$. Using our
trace form identification, this composite is the map taking a Hermitian matrix
with spectrum $\lambda$ to its diagonal entries. Then the A/G-S convexity
theorem tells us that the image is the convex hull of images of
$T$-fixed points. A Hermitian matrix is fixed under conjugation by all
diagonal matrices if and only if it itself is diagonal, which means its
entries are a permutation of $\lambda$. This proves Schur-Horn.

(Of course this presentation is historically unfair, as indicated at
the beginning -- the Schur-Horn theorem was a primary inspiration for
A/G-S convexity.)

\subsection{Stabilizer groups.}\label{stabgrps}
We mention one of the many easy properties of the definition of 
moment map. The {\dfn stabilizer algebra} of a point $m$ for the action of
a Lie algebra $\lie{k}$ is defined as the Lie subalgebra giving
vector fields that vanish at $m$. (Correspondingly, the flows they
generate leave $m$ fixed, so the stabilizer algebra is the Lie algebra of
the stabilizer group of $m$.)

\begin{Proposition*}
  Let $M$ be a Hamiltonian $K$-manifold, with moment map $\Phi$, and
  $m\in M$. Then the stabilizer algebra of $m$ is the perp of the 
  image of the differential $d\Phi$ mapping from the tangent space at $m$
  to $\lie{k}^*$.
\end{Proposition*}

(Again we're using the natural pairing between $\lie{k}$ and $\lie{k}^*$;
the perpendicular of a subspace of $\lie{k}^*$ is a subspace of $\lie{k}$.)

\begin{proof}
  A vector $\vec k$ is in the stabilizer algebra of $m$ iff the vector field
  induced on $M$ vanishes at $m$, iff the symplectic gradient of
  $\langle \Phi, \vec k \rangle$ vanishes at $m$, iff the differential
  $d\langle \Phi, \vec k \rangle$ vanishes at $m$, iff 
  $\langle d \Phi, \vec k \rangle$ vanishes at $m$, which means
  $\vec k$ is in the perp of the image of the differential $d\Phi$ at $m$.
\end{proof}

There are two interesting extremes of this.  One is that $m$'s
stabilizer group in $K$ is discrete if and only if $m$'s stabilizer
algebra in $\lie{k}$ is trivial if and only if the differential of
$\Phi$ at $m$ is onto.  On the other hand, $m$ is a $K$-fixed point
if and only if all of the vector fields induced by $\lie{k}$ 
vanish at $m$ if and only if the differential of $\Phi$ at $m$ is zero.

\begin{Corollary}
  Let $p\in \lie{k}^*$ be a boundary point of the image of a moment map 
  $\Phi: M\to \lie{k}^*$. Then each point in $\Phi^{-1}(p)$ is stabilized
  by at least a circle subgroup of $K$.
\end{Corollary}

\begin{proof}
  If $p$ is an extremal point, the image of the differential can't be
  onto -- some directions from $p$ lead outside the image of $\Phi$.
  So the perp to the differential is positive-dimensional, and
  generates a positive-dimensional subgroup of $K$, the connected
  component of the stabilizer group. And any positive-dimensional compact
  group contains a circle group.
\end{proof}

We include these proofs to show off the simple connection between the action
of the Lie group and the properties of a moment map. While A/G-S convexity
is based on Morse theory applied to the moment map, the results above are
much more pedestrian and use only basic differential geometry.

This doesn't all come for free, of course -- the important fact is that
we happen to be studying isospectral sets of Hermitians (which turn out
to be symplectic manifolds), and certain maps from them like ``take diagonals''
(which turn out to be moment maps). Most other equations one might like
to study in linear algebra are not statable as the vanishing of
a moment map. I personally take this framework as a guide to some of those 
linear algebra problems which are likely to have nice solutions.

\subsection{When $K$ is not a torus.}
There is a different, rather less pleasant, 
convexity result for noncommutative groups $K$.
The first example is $K=SO(3)$ the group of rotations of $\reals^3$,
and $M$ a coadjoint orbit. We can identify $\lie{so}(3)^*$ with $\reals^3$
and so $M$ is a sphere centered at the origin. Since this is not convex,
we know we'll have to look for a slightly subtler statement than A/G-S.

In the conjugation action of $U(n)$ on the space of Hermitian matrices
(not naturally a symplectic manifold),
we know two nice things already; 
every Hermitian matrix is diagonalizable, and if we insist
that the (real) diagonal entries then be in decreasing order, the
diagonal matrix is unique. Inside the space of diagonal Hermitian matrices,
which we can think of as the dual $\lie{t}^*$ to the Lie algebra of 
the maximal torus of $U(n)$, this picks out a certain cone $\lie{t}^*_+$
called a {\dfn positive Weyl chamber} for $U(n)$. Since $U(n)$ is the
only example we will need we won't give the general definition
of positive Weyl chamber, but
merely state that for every connected compact Lie group $K$, there is
an analogous group $T$ and polyhedral cone $\lie{t}^*_+$, such that each orbit
in $\lie{k}^*$ intersects $\lie{t}^*_+$ in a unique element.

\begin{Theorem*}
  Let $M$ be a compact connected Hamiltonian $K$-manifold, with 
  moment map $\Phi$. Then the intersection of the image of $\Phi$ 
  with the positive Weyl chamber $\lie{t}^*_+$ is a convex polytope.
\end{Theorem*}

This theorem is due to Kirwan, and is a fair bit harder than the case
$K$ commutative (the A/G-S theorem above). The new difficulties primarily
come from the points of the moment polytope lying on the boundary
of the positive Weyl chamber (in the Hermitian case: the Hermitian
matrices with repeated eigenvalues), which we call the {\dfn Weyl walls}.

We can ignore these to some extent:
define the {\dfn symplectic slice} of $M$ as the preimage under $\Phi$
of the interior of $\lie{t}^*_+$. It is then a theorem
that the symplectic slice is symplectic, and a Hamiltonian $T$-manifold,
whose image is the $K$-moment polytope minus the parts hitting the
Weyl walls. (Since it is usually noncompact we can't use it
to trivially reduce Kirwan convexity to A/G-S convexity.)

There is a slightly different statement of Kirwan convexity: instead
of intersecting the image of $\Phi$ with $\lie{t}^*_+$, we can compose
$\Phi$ with the $K$-invariant map $\lie{k}^* \to \lie{t}^*_+$ that takes 
an element to the unique point in $\lie{t}^*_+$ in its $K$-orbit.
(For $K=U(n)$ this takes a Hermitian matrix $M$ to the diagonal matrix with
$M$'s eigenvalues decreasing down the diagonal.)
Then Kirwan convexity says that the image of this composite map
is a convex polytope (the same as constructed before by intersecting).

With all this machinery we can give pleasant proofs of Horn's theorems 
(but not Horn's conjecture!) on the sum of two Hermitian matrices.

\begin{Theorem*}
  Let ${\mathcal O}_\lambda,{\mathcal O}_\mu$ be the spaces of Hermitian
  matrices with spectrum $\lambda,\mu$. Let 
  $e: {\mathcal O}_\lambda\times {\mathcal O}_\mu \to \reals^n$
  take a pair of matrices to the spectrum of their sum, listed in
  decreasing order. Then the image of $e$ is a convex polytope.
  Also, if $e(H_\lambda,H_\mu)$ is an extremal
  point of the image of $e$ and is also a strictly decreasing list,
  then $H_\lambda,H_\mu$ are simultaneously block diagonalizable.
\end{Theorem*}

\begin{proof}
  We have already explained how to identify
  ${\mathcal O}_\lambda$ and ${\mathcal O}_\mu$ with coadjoint
  orbits of $U(n)$. Therefore by property 3 of moment maps 
  listed above, their product has a Hamiltonian action of $U(n)\times U(n)$.
  (In fact it is a coadjoint orbit for this big group.)

  Consider the action of the diagonal $U(n)$, i.e. conjugating both
  Hermitian matrices by the same unitary matrix. Then by property 1
  this action is Hamiltonian, and we can compute its moment map as the
  transpose of the inclusion $\lie{u}(n) \into \lie{u}(n) \oplus \lie{u}(n)$,
  composed with the $U(n)\times U(n)$ moment map, which was just inclusion
  of the coadjoint orbit.

  The transpose of diagonal inclusion $V \to V \oplus V$ is summation
  $V^* \oplus V^* \to V^*$; so the moment map for $U(n)$'s action on
  ${\mathcal O}_\lambda \times {\mathcal O}_\mu$ just takes a pair
  of Hermitian matrices to their sum.

  But now we have convexity, as our map $e$ is just the map used
  in the alternate description of Kirwan's convexity theorem.

  If $e(H_\lambda,H_\mu)$ is a strictly decreasing list of real numbers,
  that means it's in the interior of the positive Weyl chamber,
  therefore in the image of the symplectic slice, which for  
  ${\mathcal O}_\lambda \times {\mathcal O}_\mu$ is the set of pairs 
  $(H_\lambda,H_\mu)$ whose sum is diagonal with decreasing entries
  (already a familiar set to people studying Horn's problem). 
  Then we can use some element of $U(n)$ to conjugate $(H_\lambda,H_\mu)$ 
  into the symplectic slice, on which $e$ is the moment map for the
  action of $T$.

  Since $e(H_\lambda,H_\mu)$ is on the boundary of the image of $e$,
  we can apply the corollary from subsection \ref{stabgrps}, 
  and determine that $(H_\lambda,
  H_\mu)$ is invariant under some circle in $T$ (and not just the scalar
  matrices, which fix all pairs $(H_\lambda, H_\mu$)). Being invariant
  under conjugation by a nonscalar diagonal matrix forces each of
  $H_\lambda$ and $H_\mu$ to be block diagonal.
\end{proof}

(The very careful reader will wonder why we bothered with the symplectic
slice in the above, since the proposition about stabilizers didn't 
require that the group be a torus. But just because we knew that 
$e(H_\lambda,H_\mu)$ was on the boundary of the image of $e$, we didn't
know $H_\lambda+H_\mu$ to be on the boundary of the image of $U(n)$'s
moment map -- and in fact it potentially wasn't, if $e(H_\lambda,H_\mu)$
was on a wall of the Weyl chamber.)

There is now an industry generalizing these convexity results to
larger contexts (such as ``Poisson'' actions), some of which are
definitely relevant to linear algebra; we don't pause to discuss
these, as we won't need them for Horn's problem.

\section{Symplectic quotients}

Given $M$ a Hamiltonian $K$-manifold with moment map $\Phi:M \to \lie{k}^*$,
we now have a number of theorems about the image of $\Phi$ -- put
differently, about which fibers of $\Phi$ are nonempty -- and the action
of $K$ on the fiber. Given this setup, and a point
$\mu \in \lie{k}^*$ fixed by $K$, define the {\dfn symplectic quotient} (or
{\dfn symplectic reduction}) of $M$ at
the level $\mu$ as $\Phi^{-1}(\mu)/K$.
There is a reason that these quotients are nicer to study than
the fibers themselves:

\begin{Theorem*}[Marsden-Weinstein, 1974]
  The (dense, open) smooth part of $\Phi^{-1}(\mu)/K$ inherits a canonical 
  symplectic structure. 
\end{Theorem*}

(The case most frequently studied is 
$\mu$ is a regular value of $\Phi$, since $\Phi^{-1}(\mu)$ is then a 
submanifold by the inverse mapping theorem, and the action of $K$ on it
has only finite stabilizers as discussed in subsection \ref{stabgrps}.)

In the case $K=T$, where the conjugation action is trivial, the
adjoint and coadjoint representations are also trivial. So any point
$\mu \in \lie{t}^*$ is fixed by $T$. In the case $K=U(n)$ only the
scalar Hermitian matrices are fixed by conjugation.

\begin{Example}\label{ex:hopf}
  Let $M=\complexes^n$, with the circle group $S^1 = \{\exp(i\theta)\}$
  acting by multiplication by phases. Identifying the Lie algebra of $S^1$ 
  with the reals, the moment map is $\Phi(\vec v) = \frac{1}{2} |\vec v|^2$.
  Then if $k>0$, the symplectic quotient
  $\Phi^{-1}(k)/S^1$ is $\CP^{n-1}$; this is the Hopf fibration,
  fibering $S^{2n-1}$ by circles. If $k=0$, the symplectic quotient
  is a point. If $k<0$, the symplectic quotient is empty.
\end{Example}

\begin{Example}
  Let $M={\mathcal O}_\lambda$ be the space of Hermitian matrices with
  spectrum $\lambda$, and $K=T$ the group of diagonal unitary matrices
  acting by conjugation, so $\Phi$ is the map taking a Hermitian matrix
  to its diagonal entries. Then $\Phi^{-1}(\mu)/T$ is the space of Hermitian
  matrices with eigenvalues $\lambda$, and diagonal entries $\mu$, 
  up to conjugation by diagonal unitary matrices. These spaces are
  studied in the doctoral theses \cite{Kn,G}, and the case 
  $\lambda=(1,1,0,\ldots,0)$ is described in detail in subsection
  \ref{ss:amusing}.
\end{Example}

\section{Flag manifolds}

Hopefully the previous section has convinced the reader that some of
the fundamental objects of interest in this theory are the isospectral
sets themselves. We will see now that these are not just real manifolds,
but complex manifolds, suggesting that they may be studiable using
complex algebraic geometry.

Define a {\dfn (partial) flag} in $\complexes^n$ as an increasing list 
$\mathcal V$ of 
subspaces $0 = V_0 < V_1 < \ldots < V_s = \complexes^n$, with relative
dimensions $\dim V_i / V_{i-1} = d_i$. If $s=n$ then $\mathcal V$
is called a {\dfn full flag}. If $s=2$ then $\mathcal V$ is just 
given by a
subspace of dimension $d_1$ (with the automatic $V_0$ and $V_2$).

Given a Hermitian matrix $H$, we can associate a 
partial flag ${\mathcal V}_H$ as follows: 
let $V_i$ be the sum of the eigenspaces corresponding to the $i$
smallest eigenvalues. This flag will be full if and only if $H$ has
no repeated eigenvalues.
Conversely, given a partial flag $\mathcal V$ {\em and an increasing 
list of $s$ eigenvalues $e_i$,} 
we can construct a Hermitian matrix $H_{\mathcal V}$ as the sum 
$$ H_{\mathcal V} = \sum_{i=1}^s e_i \cdot \hbox{[the projection onto  
$V_{i-1}^{\perp} \cap V_i$]}. $$
So in all, the space of Hermitian matrices with spectrum $\lambda$
is in 1:1 correspondence with a certain space of partial flags we
denote $Flags(d_1,\ldots,d_s)$
(where the $\{d_i\}$ are determined by 
$\lambda$'s repetition of eigenvalues), called a {\dfn flag manifold}.
(In the special case $s=2$, this is also called the {\dfn Grassmannian
of $d_1$-planes} and denoted $Gr_{d_1}(\complexes^n)$.) 

What have we gained?
The benefit of this view is that while Hermitian matrices only have an
evident action of the unitary group $U(n)$, the space of flags
$Flags(d_1,\ldots,d_s)$ has an action of the group $GL_n(\complexes)$
of all invertible $n\times n$ matrices -- applying a linear transformation
to a flag produces a new flag. Since $U(n)$ acts transitively on the
Hermitian matrices with a given set of eigenvalues, so too will this
larger group $GL_n(\complexes)$ act with one orbit -- given two flags,
there is a linear transformation taking one to the other. (And it may
even be taken to be unitary.)

The unitary matrices stabilizing a given Hermitian matrix -- say, a 
diagonal one with strictly decreasing eigenvalues -- are easy to understand;
they are just the diagonal unitary matrices $T$. From this we can conclude
that if $\lambda$ has no repeated eigenvalues, we can identify the space
of Hermitian matrices with spectrum $\lambda$ with $U(n)/T$. 
For the corresponding statement for $GL_n(\complexes)$ we must compute
the stabilizer of a flag -- say, the standard flag $(0<\complexes^1
<\complexes^2<\ldots<\complexes^n)$ -- which can be seen to be the
group $B$ of invertible upper triangular matrices. 
So we can identify $Flags(1,1,\ldots,1)$ with $GL_n(\complexes)/B$.

Since $GL_n(\complexes)$ is not just a real Lie group but a complex Lie
group, and $B$ a complex subgroup,\footnote{%
This means that its Lie algebra is invariant under multiplication by $i$.
By contrast, $U(n)$ is not a complex subgroup because $i$ times a
skew-Hermitian matrix is not always skew-Hermitian.}
we find out that $Flags(1,1,\ldots,1)$ is naturally a complex manifold.
As we will see in the sections to come, it is actually a ``complex
algebraic variety''.

We invite the reader to determine the corresponding statements for
the $s<n$, partial flag manifold case.

\section{Geometric invariant theory in the affine case}

Our interim goal is to describe the algebro-geometric analogue of
symplectic quotients. We first give the basics of affine algebraic
geometry, and the concept of quotient in that case; however for our linear
algebra applications we'll need to work through the additional 
complications of projective geometry.

Before getting into algebraic geometry, let us think about embeddings
of manifolds $M$ into $\reals^n$. On the one hand, if we choose
$n$ real-valued functions $\{f_i\}$ on $M$, such that any two points 
are distinguished by at least one of the functions, the map
$m \mapsto (f_1(m),f_2(m),\ldots,f_n(m))$ gives an injection\footnote{%
  Which is not good enough, if we're thinking about smooth manifolds $M$
  and want to model the smooth structure; we must also ask that the
  map be an immersion. This goes by the phrase ``separating points'',
  which we're already asking, and ``separating tangent vectors.''
  }
of $M$ into $\reals^n$. Conversely, if $M$ is a submanifold of $\reals^n$,
then the $n$ coordinate functions, restricted to $M$, give us $n$ functions
separating points.

By taking polynomials in those $n$ functions, we get an algebra of functions
on $M$, some quotient of the polynomial ring $\reals[f_1,\ldots,f_n]$.
For example, if $M$ is the parabola $y=x^2$ in $\reals^2$, and $f_1=x$
and $f_2=y$, then the algebra of functions they will generate will be
$\reals[f_1,f_2]/(f_1^2-f_2)$. (The function $f_1^2-f_2$ is not zero on the
whole plane, but it is the zero function on $M$.)

The manifolds we are actually interested in are flag manifolds
(and things we build from them), which we saw in the last section are
complex. For this reason the only rings we will bother 
to consider henceforth will be quotients of polynomial rings
with complex coefficients rather than the reals.

\subsection{The spectrum of a ring.}
Let's now think about the reverse direction: 
given a quotient $R = \complexes[f_1,\ldots,f_n]/{\mathcal I}$
of polynomial ring by an ideal, we can define {\dfn the spectrum}%
\footnote{%
If $T$ is a linear transformation $V\to V$, and $R = \complexes[T]/($
the characteristic polynomial of $T)$, then $\Spec R$ is the spectrum in
the usual sense. It is rather amusing to trace in this way
the modern-day algebraic geometry terminology back to Rydberg lines of atoms!}
$\Spec R$ as the subset of $\complexes^n$ where all the polynomials in the
ideal $\mathcal I$ vanish.

Of course, this is a bad definition, as the notation $\Spec R$ only refers
to the abstract ring $R$, and not the particular way of presenting it
as a quotient of a polynomial ring. (There may be many of these, just
as there are many embeddings of $M$ into affine spaces, of many dimensions.)
There is an alternate, equivalent, definition
of $\Spec R$ as the set of maximal ideals\footnote{%
One of the great insights of algebraic geometry in this century is that
for almost all purposes one should instead work with the prime ideals,
a somewhat larger set than the maximal ideals. But we will be avoiding 
any of the contexts in which this distinction becomes important.}
of $R$, which doesn't require us to choose an embedding.
But just as we generally deal with a manifold by picking coordinates
on it, we generally deal with {\dfn affine varieties} $\Spec R$ by
embedding them in affine space.

Given a map on rings $R \to S$, there is an induced map\footnote{For
the rings we will consider, which are finitely generated; geometrically,
this corresponds to us being able to embed our varieties in
finite-dimensional vector spaces}
from $\Spec S$ to $\Spec R$, taking a maximal ideal to its preimage.
The most fundamental example of this is $\complexes[x_1,\ldots,x_n]
\onto R$, which induces the inclusion $\Spec R \into \complexes^n$.
Not every function from $\Spec S$ to $\Spec R$ arises in this way,
only those preserving in some sense the structure of algebraic variety.
While it is possible to give a definition of such functions without
explicit reference to $R$ and $S$ we will not need this.

\subsection{Group actions.}
Now take the situation that a complex Lie group $G$ acts on our ring
$R$ by ring automorphisms. For example, if $R=\complexes[x_1,\ldots,x_n]$ 
(so $\Spec R$ is just $\complexes^n$ itself), then $GL_n(\complexes)$
acts on $R$ by linear changes of the $\{x_i\}$, and the induced action
on higher-order polynomials.

Since $G$ acts by ring automorphisms, it takes maximal ideals to maximal
ideals, so acts on the set $\Spec R$. We will be interested in forming
the quotient $(\Spec R)/G$, with the natural map $\Spec R \onto (\Spec R)/G$.

Since we're thinking in terms of algebraic geometry, this really means
we're looking for a ring $S$ with $\Spec S = (\Spec R)/G$;
equivalently, $S$ should be the ring of functions on $(\Spec R)/G$.
Each element of $S$ therefore
pulls back to a $G$-invariant function on $\Spec R$. 
So $S$ will have to map back to the $G$-invariant functions on $\Spec R$,
which we denote $R^G$; this is the clear candidate for $S$.
Since $R^G$ is the ``ring of invariants'', 
we will call $\Spec (R^G)$ the {\dfn geometric invariant theory quotient}
(or GIT quotient) of $\Spec R$ by $G$, and denote it $(\Spec R) // G$.

It is not quite true, though, that $(\Spec R)//G = (\Spec R)/G$; rather we only
have a map $(\Spec R)/G \to (\Spec R)//G$. It is already a tricky theorem that
when $G$ is ``reductive'' (like $GL_n(\complexes)$, and all the other
groups we will be using) that (1) $R^G$ is finitely generated, 
so $(\Spec R)//G$ can be embedded in a finite dimensional space, and (2)
the natural map from $(\Spec R)/G \to (\Spec R)//G$ is onto.

\begin{Example}\label{ex:rescale}
  Let $G=\complexes^\times$, the nonzero complex numbers under multiplication,
  acting on $R=\complexes[x_1,\ldots,x_n]$ by rescaling each coordinate
  the same way. Then the induced action on $G$ on $\complexes^n$ is
  also by rescaling. The ordinary set-theoretic quotient is not Hausdorff;
  it's projective space union the point $\vec 0$ in the closure of every
  other point. The only invariant polynomials in $R^G$ are the constants,
  so $\Spec (R^G) = pt$. (This will turn out to be related to the
  $k=0$ case of example \ref{ex:hopf}.)
\end{Example}

Obviously this is unsatisfying; we'd rather the quotient be projective space
itself (somehow including the $k>0$ case of example \ref{ex:hopf}). 
But projective space is not $\Spec S$ for any $S$, by Liouville's
theorem -- any function on projective space is constant! We will need 
subtler constructions than $\Spec$ to make projective spaces, and more
to the point, flag manifolds.

\section{Geometric invariant theory in the projective case, and
  the Kirwan/Ness theorem}

Let $R =\oplus_{k\in \naturals} R_k$ now be a {\em graded} ring, meaning
that the product of an element of $R_k$ with an element of $R_m$ lands in
$R_{k+m}$. The standard example is $R$ a polynomial ring, with $R_k$
the homogeneous polynomials of degree $k$.

We will define $\Proj R$ in several equivalent ways. The simplest,
least useful for visualizing examples, is that $\Proj R$
the set of maximal {\em graded} ideals of $R$, where a graded ideal
is one equal to the direct sum of its intersections with the $R_n$.

For the second, note that $R$ has a natural action of $\complexes^\times$,
acting on $R_k$ by rotating it with speed $k$. There is a natural 
$\complexes^\times$-invariant surjection $R \onto R_0$ (which would not
be true if $R$ were $\integers$-graded instead of $\naturals$-graded).
So there is a map backwards $\Spec R_0 \into \Spec R$.
Then 
$$ \Proj R := (\Spec R \smallsetminus \Spec R_0) / \complexes^\times.$$

\begin{Example}
  Let $R=\complexes[x_1,\ldots,x_n]$ with $R_n$ the homogeneous polynomials
  of degree $n$. Then as discussed in example \ref{ex:rescale}, $R_0$
  is just the constants, and the point $\Spec R_0$ includes into 
  $\Spec R = \complexes^n$ as the origin $\vec 0$. Then 
  $\Proj R = (\complexes^n \smallsetminus \{\vec 0\})/\complexes^\times$
  is just the usual definition of $\CP^{n-1}$.
\end{Example}

So $\Proj$ of a polynomial ring (with each variable given degree 1) 
is projective space, just like $\Spec$ of a polynomial ring is affine space.
And just as writing a ring $R$ as the quotient of a polynomial ring,
$\complexes[x_1,\ldots,x_n] \onto R$, dually gives us an inclusion of
$\Spec R \into \complexes^n$, writing a graded ring $R$ as the quotient
of a polynomial ring by a graded ideal $\mathcal I$ dually gives us an
inclusion $\Proj R \into \CP^{n-1}$. This gives a third description
of $\Proj R$, 
in the case $R$ is presented as $\complexes[x_1,\ldots,x_n]/{\mathcal I}$
and all the $x_i$ are degree 1;
as the points in $\CP^{n-1}$ where the elements of ${\mathcal I}$ all vanish.
In this case $\Proj R$ is called a {\dfn projective variety.}

The description we will use most is the second, in terms of $\Spec R$
and $\Spec R_0$.

\subsection{Another relation between $\Spec$ and $\Proj$.} Given an
ungraded ring $R_0$, we can define $R := R_0[l]$, where $l$ is a new variable
assigned formal degree 1. Then $\Spec R = \Spec R_0 \times \complexes$, so
$\Spec R \smallsetminus \Spec R_0 = \Spec R_0 \times \complexes^\times$,
and $\Proj R = \Spec R_0$. Upshot: anything we can make as a $\Spec$,
we can also make as a $\Proj$.

\subsection{Maps of graded rings.} Given a ring homomorphism $R \to S$,
we defined a natural map from $\Spec S \to \Spec R$. So it is natural
to assume that a grading-preserving map $R \to S$ induces a map on
the corresponding $\Proj$'s. This is not the case! Rather:

\begin{Proposition}
  Let $f: R\to S$ be a homomorphism of graded rings. Let $X$ be the
  points of $\Spec S\smallsetminus \Spec S_0$ such that every function in
  $f(R)$ vanishes at $X$, and define $(\Proj S)^{us}_f$ as the
  quotient of $X$ by the natural action of $\complexes^\times$.
  Then $f$ does not necessarily induce a map from $\Proj S \to \Proj R$,
  but does induce one from $\Proj S \smallsetminus (\Proj S)^{us}_f \to \Proj R$.  
\end{Proposition}

\begin{proof}
  To see what the problem is, let's trace through the definition of $\Proj$.
  Since $R$ and $S$ are graded, i.e. have $\complexes^\times$-actions,
  $\Spec R$ and $\Spec S$ have $\complexes^\times$ actions. The
  map between the rings being grading-preserving is equivalent to it
  being $\complexes^\times$-equivariant. Consequently, the map backwards
  on $\Spec S\to \Spec R$ is $\complexes^\times$-equivariant. 
  
  The problem comes when we try to rip out $\Spec S_0$ and $\Spec R_0$.
  Of course there is no problem when in restricting the map to
  $\Spec S \smallsetminus \Spec S_0$. But the image of this may not land
  inside $\Spec R \smallsetminus \Spec R_0$; there may be points of 
  $\Spec S \smallsetminus\Spec S_0$ that hit $\Spec R_0$. Once we rip them out,
  then there's no problem.
\end{proof}

In the context to come, the set $(\Proj S)^{us}_f$ will be called
the {\dfn unstable set} (hence the $us$).

\subsection{Geometric invariant theory.}
We now apply what we've learned about maps of $\Proj$'s to the case
we studied in the last section, $R^G \into R$. Given a graded ring $R$
with an action of a group $G$, define the {\dfn geometric invariant
theory quotient}
$$ (\Proj R) // G := \Proj (R^G). $$
There is a natural map from $\Proj R$ minus the unstable set
described by the proposition above (the set of points where all
$G$-invariant functions vanish) to $(\Proj R)//G$.

While this may seem like a totally canonical definition, there are two
traps for the unwary. In many problems one starts with the variety
$\Proj R$ and its $G$-action rather than the ring $R$, and does not
want the extra choice that comes in finding an $R$ (of which, unlike
in the affine case, there
will be many). And even when one has chosen an $R$, there may be many
ways to get $G$ to act on it inducing the same action on $\Proj R$.
It turns out that both these choices matter -- which is to say,
the notation $(\Proj R) // G$ is misleading. 

\begin{Example}\label{ex:CPnGIT}
  Let $X=\complexes^n$, and $G=\complexes^\times$ act by rescaling.
  Then $X=\Spec \complexes[x_1,\ldots,x_n]
  =\Proj \complexes[x_1,\ldots x_n,l]$, where the $\{x_i\}$ are all
  degree $0$, and $l$ is degree $1$. The action of $G$ on our homogeneous
  coordinate ring $R$ is not quite determined -- while the $x_i$ are required
  to all be weight $1$, the action on the variable $l$ can be any weight $k$,
  i.e. $z \cdot l = z^k l$.

  If $k>0$, there are no invariant polynomials other than constants: 
  $R^G = \complexes$, and $\Proj (R^G)$ is the empty set.

  If $k=0$, the invariant ring $R^G = \complexes[l]$, and $\Proj (R^G)$ 
  is a point.

  If $k=-1$, the invariant ring $R^G=\complexes[x_1 l,x_2 l,\ldots, x_n l]$
  with all products of degree 1, so $\Proj (R^G) = \CP^{n-1}$. It is
  a little trickier to see that this is true for {\em all} negative $k$.
  (It is interesting to compare this to example \ref{ex:hopf}. The
  difference in signs can be explained but we do not do so here.)
\end{Example}

\begin{Example}\label{ex:GrGIT}
  Let $X=\complexes^{m\times n}$, the space of $m\times n$ matrices,
  with $m\leq n$, and let $G=GL_m(\complexes)$ act on $X$ on by left
  multiplication.  Then as in the previous example (the $m=1$ case of
  this one) there is a parameter $k$ describing the action of $G$ on
  the extra variable $l$, such that for $k>0$ the quotient $X // G$ is
  empty, and for $k=0$ the quotient is a point.

  Regard $X$ as the space of $m$ vectors in $\complexes^n$.
  For $k<0$, the unstable set turns out to be 
  the set of linearly dependent $m$-tuples,
  and the GIT quotient $X//G$ is thus the Grassmannian $Gr_m(\complexes^n)$. 
  (If you take a linearly independent set
  of $m$ vectors, and quotient by the action of $GL_m(\complexes)$, 
  you forget the vectors and only remember the $m$-subspace they span.)
\end{Example}

\subsection{The equivalence of symplectic and GIT quotients.}
We spell out a moment map that's a special case of things already said:
complex projective space $\CP^{n-1}$ is a Hamiltonian $U(n)$-manifold,
with moment map
\begin{align*}
\CP^{n-1} &\to \lie{u}(n)^* \\
[\vec v] &\mapsto [\hbox{the rank $1$ projection onto $\complexes\vec v$}]
\end{align*}
using the identification already given between $\lie{u}(n)^*$ and
the space of Hermitian matrices. In fact we are observing here that
$\CP^{n-1}$ is a coadjoint orbit of $U(n)$, 
the orbit ${\mathcal O}_{(1,0,\ldots,0)}$.

It is an easy fact,\footnote{%
once one proves the non-easy Darboux theorem, given the unusual way we've
defined symplectic} 
not proven here, that any complex submanifold of
$\CP^{n-1}$ inherits a symplectic structure from $\CP^{n-1}$.
If $X \subseteq \CP^{n-1}$ is preserved by the action of 
a subgroup $K\leq U(n)$, we can compute the moment map using the
properties we listed of moment maps:
$$ X \into \CP^{n-1} \to \lie{u}(n)^* \to \lie{k}^*. $$

Lastly, given an action of a compact group $K$ on a complex vector space, 
there is a unique extension to an action of {\dfn the complexification
$K^\complexes$ of $K$}. We will not stop to define this group in detail;
suffice it to say that every compact group is a subgroup of a unique
complex group, such that the Lie algebra of the complex group is the
complexification of the Lie algebra of the compact group. The only
example that will interest us is $U(n)^\complexes = GL_n(\complexes)$
(every matrix is uniquely the sum of a skew-Hermitian and $i$ times
a skew-Hermitian).

We are ready to state the deepest theorem in this paper,
proven separately by Kirwan and Ness:

\begin{Theorem}
  Let $K$ act on the graded ring $R = \complexes[x_1,\ldots,x_n]/{\mathcal I}$,
  all generators $x_i$ of degree $1$.
  So $\Proj R$ is a subvariety of $\CP^{n-1}$, and (if smooth)
  a Hamiltonian $K$-manifold. Let $\Phi$ be the moment map calculated above
  for $K$'s action on $M$. Then there is a natural identification
  $$ \Phi^{-1}(0)/K \quad\iso\quad (\Proj R)//K^\complexes $$  
  between the symplectic quotient and geometric invariant theory quotient.
\end{Theorem}

In fact there is a more general result, slightly more complicated to
state, in which $\Proj R$ is not necessarily projective.
We already saw this in the case of $K=S^1$ acting on $\complexes^n$,
in examples \ref{ex:hopf} and \ref{ex:CPnGIT}.

\section{Back to flag manifolds}

Since we've only stated the symplectic quotient equals GIT quotient 
theorem at the level $\Phi=0$, we need to slightly modify the Horn
problem $A+B=C$ to $A+B+C=0$. 
Let $M={\mathcal O}_\lambda \times {\mathcal O}_\mu \times {\mathcal O}_\nu$, 
the space of triples of Hermitian matrices with eigenvalues
respectively $\lambda,\mu,\nu$. Then the moment map for the diagonal
action of $U(n)$ on the product is the sum of the three matrices, and
Horn's problem asks when this symplectic quotient is nonempty.

To jack into the symplectic vs. GIT theorem, we need to express 
${\mathcal O}_\lambda$ as $\Proj$ of something, or equivalently explain
how to embed it in projective space such that the restriction of the
symplectic structure on projective space matches that on 
${\mathcal O}_\lambda$. This turns out to only be possible if
$\lambda$ is integral, and involves some very classical geometry.

\subsection{The $U(n)$ case of the Borel-Weil-Bott-Kostant theorem.}
In fact we will not need the details presented in this section
in what follows, but it lets us avoid some
representation theory of $GL_n(\complexes)$.

Recall the flag manifolds $Flags(d_1,\ldots,d_s)$ defined before; we
will restrict to the case $d_i \equiv 1$ and just write
$Flags(\complexes^n)$, leaving the interested reader to work out the
general case. Denote by $Gr_k(\complexes^n)$ the Grassmannian of
$k$-dimensional subspaces of $\complexes^n$ (these are very partial
flag manifolds). There is a natural forgetful map from $Flags(\complexes^n)$
to each Grassmannian, taking a flag to its subspace of dimension $k$,
forgetting all the subspaces below and above.

Since a flag is just a list of subspaces, the product map
$$ Flags(\complexes^n) \to \prod_{k=1}^{n} Gr_k(\complexes^n) $$
is an inclusion.

Given a $k$-dimensional subspace $A$ of $\complexes^n$, we can wedge
together a basis $\{\vec a_i\}$ 
of $A$ to get a nonzero alternating tensor in $\wedge^k \complexes^n$.
If we change the basis, our element $\wedge_{i=1}^k \vec a_i$ only
changes by a scalar factor, so gives a well-defined element of the
projective space. This is called the {\dfn Pl\"ucker embedding} 
$$ Gr_k(\complexes^n) \to \PP(\wedge^k \complexes^n). $$

Given a vector $\vec v \in V$, and a natural number $a$, 
we can tensor $\vec v$ with itself $a$ times to get a symmetric tensor
$\vec v^{\tensor a} \in \Sym^a(V)$. This descends to the
projectivized spaces $\PP V \into \PP(\Sym^a V)$, 
and is called the {\dfn $a$th Veronese embedding}. 
Choosing a family $\{a_k\}$ of naturals, we get maps
$$ \PP(\wedge^k \complexes^n) \to \PP(\Sym^{a_k} (\wedge^k \complexes^n)).$$
In the very special case that $\dim V=1$, we don't need $a$ to be a
natural number -- for $a<0$ we define $\Sym^a(V) := \Sym^{-a}(V^*)$.
This will be handy in the case $k=n$, where $\wedge^n \complexes^n$ is
one-dimensional.

Finally, given vectors $\vec v\in V$ and $\vec w\in W$, we can tensor
them together, inducing a map on projective spaces 
$\PP V \times \PP W \to \PP(V\tensor W)$ called the {\dfn Segre embedding}.  
In the case at hand, this gives a map
$$ \prod_{k=1}^{n} 
\PP(\Sym^{a_k} (\wedge^k \complexes^n))
\to
\PP(\tensor_{k=1}^{n} \Sym^{a_k} (\wedge^k \complexes^n)).$$

\begin{Proposition*}[one small aspect of the Borel-Weil-Bott-Kostant theorem]
  Let $\lambda = (\lambda_1> \ldots >  \lambda_n)$ be a strictly
  decreasing list of integers, and $a_k = \lambda_k - \lambda_{k+1}$ for
  $k=1\ldots n$ (taking $\lambda_{n+1}$ to be zero). 
  Then the composite of the above maps
$$ Flags(\complexes^n) 
\to
\PP(\tensor_{k=1}^{n} \Sym^{a_k} (\wedge^k \complexes^n))$$
induces a symplectic structure on $Flags(\complexes^n)$ matching the one 
given by its diffeomorphism with the coadjoint orbit ${\mathcal O}_\lambda$.
\end{Proposition*}

(Again, we exhort the reader to think about the case of $\lambda$
only weakly decreasing, the partial flag manifold case.)

Since we have exhibited ${\mathcal O}_\lambda$ as a variety embedded
in projective space, we can look at its homogeneous coordinate ring,
which we will denote $R^\lambda$, reserving subscripts to indicate
graded pieces.

\subsection{The Borel-Weil theorem.}
It is evident that all the maps described above are equivariant with 
respect to the action of $GL_n(\complexes)$ on $\complexes^n$, 
so $GL_n(\complexes)$ therefore acts on the homogeneous coordinate rings
$R^\lambda$, and in particular each graded piece $(R^\lambda)_k$ is
a representation of $GL_n(\complexes)$. We pause to mention three important
facts about these representations.

\begin{Proposition}[part of Borel-Weil]
  \begin{enumerate}
  \item Each graded piece $(R^\lambda)_k$ of each homogeneous coordinate
    ring for the flag manifold is an irreducible representation of 
    $GL_n(\complexes)$.
  \item $(R^\lambda)_k$ is isomorphic, as a representation, to
    $(R^{k\lambda})_1$.
  \item Every irreducible representation of $GL_n(\complexes)$ arises as
    $(R^\lambda)_1$ for a unique $\lambda$.
  \end{enumerate}
\end{Proposition}

We denote $(R^\lambda)_1$ by $V_\lambda$; for reasons we will not go
into here it is called the {\dfn irreducible representation of
$GL_n(\complexes)$ of highest weight $\lambda$}. So we can identify
$R^\lambda$ with $\oplus_{n\in \naturals} V_{n\lambda}$.

A little bit more can be said, for those who know the classification of
irreducible representations of Lie groups in terms of ``highest weights'',
but as representation theory is not in the title of this paper we will
not explore this further.

\section{Easy consequences}

We need one more fact about the homogeneous coordinate rings 
$R^\lambda$ for the flag manifold: they have no zero divisors.
This is implied, though we will not explore how, by the 
(much stronger) fact that the flag manifold is connected and smooth.

We're now ready to apply our big gun, the equivalence of symplectic
and GIT quotients.

\begin{Theorem*}
  Let $\lambda,\mu,\nu$ be weakly decreasing sequences of integers.
  Then the space 
  $$ \bigg\{ (H_\lambda,H_\mu,H_\nu): H_\lambda+H_\mu+H_\nu = 0\bigg\} \bigg/ U(n)$$
  can be identified with
  $$ \Proj \bigoplus_{k\in \naturals} (V_{k\lambda} \tensor V_{k\mu} \tensor
  V_{k \nu})^{GL_n(\complexes)}. $$
\end{Theorem*}

\begin{proof}
  We've already seen that the first space is the symplectic quotient of
  ${\mathcal O}_\lambda\times{\mathcal O}_\mu\times{\mathcal O}_\nu$ 
  by the diagonal action of $U(n)$, at the level $0$. 
  By the symplectic vs. GIT equivalence, this is the geometric invariant
  theory quotient of a certain product of flag manifolds by 
  $GL_n(\complexes)$.

  The Borel-Weil-Bott-Kostant theorem, or rather the small part of it
  presented here, explains how to embed 
  ${\mathcal O}_\lambda,{\mathcal O}_\mu,{\mathcal O}_\nu$ into projective
  space in order to be able to apply the symplectic vs. GIT theorem.
  The individual coordinate rings are $\oplus_{n\in \naturals} V_{n\lambda}$
  (likewise $\mu,\nu$); the coordinate ring of the product is, by 
  Segre embedding, 
  $ \oplus_{n\in \naturals} (V_{n\lambda} \tensor V_{n\mu} \tensor
  V_{n\nu})^{GL_n(\complexes)}. $

  Then to take the GIT quotient we take $\Proj$ of the invariant subring.
\end{proof}

This explains a couple of connections already noticed in the
literature between the Horn's problem and the representation theory of
$GL_n(\complexes)$:

\begin{Corollary}
  Let $V_\lambda,V_\mu,V_\nu$ be three irreducible representations of
  $GL_n(\complexes)$. If there is a $GL_n(\complexes)$-invariant vector
  in the triple tensor product, then there exist Hermitian matrices
  $H_\lambda, H_\mu,H_\nu$ with the corresponding spectra and zero sum.
\end{Corollary}

\begin{proof}
  If the tensor product has an invariant vector, then the coordinate
  ring of the GIT quotient is nontrivial in its degree 1 piece.
  Since this ring is a subring of a ring with no zero divisors,
  it itself has no zero divisors, and therefore is nontrivial in
  all degrees. This is enough to conclude that $\Proj$ of it, 
  the GIT quotient, is nonempty.
  
  By the identification in the theorem, the symplectic quotient is
  therefore nonempty. So there are three Hermitian matrices with the
  desired spectra adding to zero.
\end{proof}

\begin{Corollary}
  Let $\lambda,\mu,\nu$ be weakly decreasing lists of integers.  
  If there exist Hermitian matrices $H_\lambda, H_\mu,H_\nu$ with the
  corresponding spectra and zero sum, then for some $k>0$,
  the tensor product $V_{k\lambda} \tensor V_{k\mu} \tensor
  V_{k\nu}$ has a $GL_n(\complexes)$-invariant vector.
\end{Corollary}

\begin{proof}
  The condition says the symplectic quotient is nonempty, so running
  the theorem in reverse, the GIT quotient is nonempty.
  Consequently the ring of invariants is nontrivial. Therefore in some
  graded piece $n>0$ it is nontrivial, giving the desired result.
\end{proof}

(Both of these corollaries appear in Klyachko's work; he essentially
repeats the proof of Kirwan/Ness in this special case.)

Obviously there is a mismatch here; we would like to know that the
Hermitian problem has a solution if and only if the tensor product
has an invariant vector. But the big machine presented here is not
powerful enough to rid us of this $k$. Other techniques are necessary
and the first proof that $k$ can indeed be taken to be $1$ appeared
in \cite{KT}.

\section{Conclusions, and an amusing example}

We saw that certain famous equations in linear algebra, such as
$H_\lambda + H_\mu = H_\nu$ (or better, 
$H_\lambda + H_\mu + H_{-\nu} = 0$), are the conditions that a
certain moment map be equal to zero. Whenever this happens we can plug
into the theory of Hamiltonian actions on symplectic manifolds.
The most spectacular results here establish the convexity properties
of the images of these moment maps, and go some way toward 
determining the image. 

In fact the details of the proofs are not actually very different from
the hands-on techniques used e.g.  by Horn himself; the benefit here
is in establishing a framework that points out which problems are
likely to accede to such local analysis, and in particular which are
likely to lead to linear inequalities and polytopes.

In the (frequent) case that the symplectic manifolds under study are
algebraic varieties, one then has a totally new viewpoint, replacing
the spaces by homogeneous coordinate rings, and the study of the
moment map by invariant theory. This connects up the symplectic problem
with representation theory problems, illuminating both.

Lastly, we emphasize another lesson (not really applied in this paper)
from both symplectic and algebraic geometry; the questions like ``for
which $\lambda,\mu,\nu$ is the following symplectic quotient
nonempty'' should really be viewed as just the first step in a more
detailed study of the symplectic quotient itself. What is its
dimension? symplectic volume? its Betti numbers, and cohomology ring?

\subsection{An amusing example.}\label{ss:amusing}
Let $X=\complexes^{2\times n}$, the space of $2\times n$ matrices, 
with its left action of $U(2)$ and right action of $U(n)$.
These actions are Hamiltonian with moment maps $M \mapsto -M M^*$ and
$M \mapsto M^* M$. We will only be interested in the action of the 
diagonal subgroup $T^n \leq U(n)$, whose moment map picks out the
diagonal entries of $M^* M$. Note that $\Tr M M^* = \Tr M^* M$, so
the two moment maps are not entirely unrelated.

We will symplectic quotient by both groups $U(2)$ and $T^n$.
The symplectic quotient by $U(2)$ at level $s {\bf 1}$ asks that our
$2$ vectors in $\complexes^n$ be orthogonal with
norm-square $s$. So for $s>0$, the symplectic quotient is the
Grassmannian of $2$-planes in $\complexes^n$. (Which is good,
because we got this answer when we did this earlier by GIT,
in example \ref{ex:GrGIT}.)
Since that is a coadjoint orbit of $U(n)$, studying its reductions
by $T^n$ is really the Schur-Horn problem we first mentioned.

If we consider this in the opposite order, we get a very different
picture. In symplectic-quotienting $\complexes^{2\times n}$ by $T^n$,
each diagonal entry in $T^n$ acts on its own copy of $\complexes^2$.
The symplectic quotient of $\complexes^2$ by $U(1)$ 
is either empty, a point, or the Riemann sphere 
$\CP^1$, depending on the level set $\vec a = (a_1,\ldots,a_n)$ chosen
in $\lie{t}^n \,\iso\, \reals^n$. Let's take each $a_i>0$ so the 
symplectic quotient by $T^n$ is a product of $\CP^1$'s.
In fact the moment map for the residual action of $U(2)$ on each of
these $\CP^1$ identifies it with a sphere $S^2_{a_i}$ of radius $a_i$.

Think of such a sphere as the set of steps in $\reals^3$ of length $a_i$;
then the product $\prod_{i=1}^n S^2_{a_i}$ can be thought of as the 
set of $n$-step polygonal paths in $\reals^3$, with $i$th step length $a_i$.
Since the moment map for the diagonal action of $U(2)$ on this product
is the sum of the individual moment maps, we can think of it as the
function taking a path to its endpoint.

The symplectic quotient at level zero is then the set of polygons in
$\reals^3$ (because the moment map condition requires that the
path terminate at the origin), with edge lengths $\{a_i\}$, 
considered up to rotation. The connection of polygon spaces 
to the Schur-Horn problem was first noted in \cite{HK}.

\bibliographystyle{alpha}    

\end{document}